\documentclass[12pt,aps,onecolumn]{revtex4}

\usepackage{amsmath,amssymb}

\begin{document}

\title{\bf
The Eratosthenes Progression \\
$p_{k+1}=\pi^{-1}(p_k), k=0,1,2,...,~p_0 \in N \setminus P$\\
Determines An Inner Prime Number Distribution Law
\cite{MTCP}}

\author{Lubomir Alexandrov}

\email{alexandr@thsun1.jinr.ru}

\affiliation{
Bogoliubov Laboratory of Theoretical Physics,
{\rm Joint Institute for Nuclear Research, 141980 Dubna, Russia}}


\begin{abstract}
\textbf{Abstract}:
Formulas of $\pi(x)$-fine structure are presented.
\end{abstract}

\maketitle

Let $\pi(p_n)$ denotes $\# n$ of prime
$p_n \in P=\{2,3,5,7,11,13...\}~~(\pi(p_n)=n)$ and $\pi^{-1}(n)$ denotes the
prime $p_n$ with $\# n~~
(\,\, \pi^{-1}(n)=p_n; \,\,
\pi^{-1}(n) \equiv \mbox{Prime}[n], \, Mma) $.

For the sets of primes (``Eratosthenes rays'')

\vspace{-3mm}

\begin{equation} \label{eq1}
r_{p_0} = \{ p_{k+1} = \pi^{-1} ( p_k ): p_0 \in N,  k=0,1,2,... \}
\end{equation}

\vspace{-3mm}

the assertions

\vspace{-3mm}

\begin{equation} \label{eq2}
\bigcap_{p_0 \in \overline{C}} r_{p_0} = \O, ~~~
\bigcup_{p_0 \in P} r_{p_0} \subset \bigcup_{p_0 \in \overline{C}} r_{p_0},~~~
\bigcup_{p_0 \in N} r_{p_0} \subset \bigcup_{p_0 \in \overline{C}} r_{p_0},~~~
P = \bigcup_{p_0 \in \overline{C}} r_{p_0},
\end{equation}
(\, where $\overline{C} = N \setminus P=\{ 1 \} \bigcup C,~ N$ is a set of
naturals and  $C$ is a set of composites \,) are~true
\cite{Sofia1964,NonAs}.

From (\ref{eq2}) we obtain matrix representations of the prime and natural numbers:
\begin{center}
$^2P = \{ r_{p_0} \}_{p_0 \in \overline{C}} = \{ p_{\mu \nu} \}~~~~~~ {\rm
and}~~~~~~^2N=\{\overline{C}, \, \,^2P \}.$
\end{center}

\newpage

The left upper corner of matrix $^2N$ looks like \cite{NonAs}, p. 18 :

\bigskip

$^2 N = \left [
\begin{array}{ccccccccccc}
       1 & 2      & 3  & 5   & 11   & 31     & 127    & 709 & 5381 & 52771 \ldots \\
       4 & 7      & 17 & 59  & 277  & 1787   & 15299  & \ldots  & &  \\
       6 & 13     & 41 & 179 & 1063 & 8527   & \dots  & & &  \\
       8 & 19     & 67 & 331 & 2221 & 19577  & \ldots & & &   \\
  \vdots & \vdots & \vdots & \vdots & \vdots & \vdots & & &
\end{array}\right ].$

\bigskip

Let $\pi(n',n''),~ n', n'' \in N$ denotes the number of primes within a given
interval $(n', n'')$.
The elements of matrices $^2P$ and $^2N$ satisfy the relations:
\bigskip

$~~p \in ~^2P_{-1} \Longleftrightarrow \pi(p) \in P,
~~~~~ p \in P_1 \Longleftrightarrow \pi(p) \in \overline{C}$,

~~~~~~~~~~~~~~~where $^2P_{-1} =~^2P~\setminus P_1,~ P_1 = {\rm column}({p_{\mu 1}})$;
\vspace{-6mm}

\begin{equation}\label{eq3}
\left .
\begin{array}{l}
\hspace{-8mm}\pi(p_{\mu 1},0) = p_{\mu 0} -1, \quad \mu \geq 1, \, p_{\mu 0} \in \overline C ;  \\
\hspace{-8mm}\pi(p_{\mu \nu_1},p_{\mu \nu_2})=p_{\mu, \nu_2 -1} - p_{\mu, \nu_1 -1} - 1,~~ \mu, \nu_1 \geq 1,~ \nu_2=\nu_1+\alpha,~ \alpha=1,2,3...; \\
\hspace{-8mm}\pi(p_{\mu_1 \nu_1},p_{\mu_2 \nu_2})=| p_{\mu_1, \nu_1 -1} - p_{\mu_2, \nu_2 -1}| - 1,~~ \mu_i, \nu_i \geq 1,~ i=1,2,~\mu_1\neq \mu_2.
\end{array}
\right \}
\end{equation}

\bigskip

The differences
$p_{p_0,k+1}-p_{p_0,k} > p_{p_0,k}(\ln{p_{p_0,k}}-1),~~ k=0,1,2,...,~ p_0=2, p_0 \in C$
are monotonically increasing along the raws of $^2P$ and
the explicit law (\ref{eq1}) together with its
corollaries (\ref{eq3}) of the $\pi(x)$--fine structure is valid.
From there a possibility follows to construct
the prime number spiral and web on the plane
$R^2$ re-creating the prime number distribution in details.

Applications of these results to the identification problem
of new transactinides as well as to some actual quantum physics
and molecular biology problems are suggested.

\end{document}